\documentclass{amsart}
\usepackage{amssymb, amsmath, latexsym, mathabx, dsfont, enumerate}

\renewcommand{\baselinestretch}{\baselinestretch}
\renewcommand{\baselinestretch}{1.1}
\numberwithin{equation}{section}

\newtheorem{thm}{Theorem}[section]

\newtheorem{rmk}[thm]{Remark}

\newcommand{\Mod}[1]{\ (\mathrm{mod}\ #1)}

\begin{document}

\title{Ternary quadratic forms representing same integers}
\author{Jangwon Ju}

\address{Department of Mathematics, University of Ulsan, Ulsan 44610, Republic of Korea}
\email{jangwonju@ulsan.ac.kr}
\thanks{This work was supported by the National Research Foundation of Korea(NRF) grant
funded by the Korea government(MSIT) (NRF-2019R1F1A1064037).}



\begin{abstract}
 In 1997, Kaplansky conjectured that if two positive definite ternary quadratic forms with integer coefficients have perfectly identical integral representations, then they are isometric, both regular, or included either of two families of ternary quadratic forms. 
In this article, we prove the existence of pairs of ternary quadratic forms representing same integers which are not in the Kaplansky's list.
\end{abstract}

\maketitle

\section{Introduction}

For a positive definite quadratic form $f$ of rank $r$ with integer coefficients and a nonnegative integer $n$, we define a set 
$$
R(n,f)=\{(x_1,x_2,\dots,x_r)\in\mathbb{Z}^r \mid f(x_1,x_2,\dots,x_r)=n\}.
$$ 
We also define $r(n,f)=|R(n,f)|$. 
It is well-known that if $f$ is positive definite, then $r(n,f)$ is always finite. 
Determining $r(n,f)$ for an arbitrary positive definite quadratic form $f$ is  quite an old problem which is still widely open. 

For a quadratic form $f$ of rank $r$, the theta series of $f$ is given by 
$$
\theta_f(z)=\sum_{\mathbf{x}\in\mathbb{Z}^r} q^{f(\mathbf{x})}=\sum_{n=0}^{\infty}r(n,f)q^n\quad (q=e^{2\pi i z}).
$$
It is well known that $\theta_f(z)$ is a modular form of weight $\frac r2$.
If two quadratic forms $f$ and $g$ are isometric over $\mathbb{Z}$, then clearly $\theta_f(z)=\theta_g(z)$. 
However, the converse is not true in general.
In \cite{wit}, Witt proved that there are two positive definite even unimodular quadratic forms of dimension 16
which have the same theta series.
After that, such examples were found even in dimension 12, 8, and 4 (for details, see \cite{kne}, \cite{ki2} and \cite{sch1}).
In \cite{cs}, Conway and Sloane found a 4–parameter family of pairs of quaternary quadratic forms with same theta series. 
Furthermore, it was proved that the Conway-Sloane tetralattice pairs are not isometric over $\mathbb{Z}$ (for details, see \cite{ch}).
However, it is well known that unary and binary positive definite quadratic forms are determined up to integral equivalence by their theta series. 
Furthermore, Schiemann proved that if two positive definite ternary quadratic forms have same theta series, then they are isometric over $\mathbb{Z}$ (for details, see \cite{sch2}).

For a positive definite quadratic form $f$ with integer coefficients, we define the set $Q(f)$ of all nonnegative integers $n$ such that $f(x,y,z)=n$ has an integer solution. 
For quadratic forms $f_1,~f_2,~g_1$ and $g_2$, we write 
$
\{f_1,~g_1\}\simeq\{f_2,~g_2\}
$ 
if $f_1\simeq f_2$ and $g_1\simeq g_2$ or $f_1\simeq g_2$ and $g_1\simeq f_2$.

In 1938, Delone proved in \cite{del} that for two binary quadratic forms $f$ and $g$, if $Q(f)=Q(g)$, then they are isometric over $\mathbb{Z}$ except the case when 
$$
\{f,~g\} \simeq \{x^2+xy+y^2,~x^2+3y^2\}.
$$
In 1997, Kaplansky conjectured in his letter to Schiemann that for two ternary quadratic forms $f$ and $g$, if $Q(f)=Q(g)$, then one of the following holds:
\begin{enumerate}[(i)]
\item $f\simeq g$;
\item $f$ and $g$ are both regular;
\item $\{f,~g\}\simeq\{ax^2+by^2+bz^2+byz,~ax^2+by^2+3bz^2\}$ for some $a,b\in\mathbb{Z}$;
\item $\{f,~g\}\simeq\{ax^2+ay^2+az^2+byz+bxz+bxy,~ax^2+(2a-b)y^2+(2a+b)z^2+2bxz\}$ for some $a,b\in\mathbb{Z}$. 
\end{enumerate}
Here, a quadratic form $f$ is called regular if it represents all integers that are locally represented by $f$. 
It was proved that the conjecture holds if only diagonal ternary quadratic forms are considered (for details, see \cite{tan}). 

Recently, Oishi-Tomiyasu gave a list of candidates of sets of ternary quadratic forms violating Kaplansky conjecture. Actually, it was proved that ternary quadratic forms contained in each set $S_i$ for $i=1,\cdots,15$ in Table 1.1 which are not isometric each other and are not regular represent same integers up to $3,000,000$. 
Furthermore, she provided a modified version of  Kaplansky conjecture (for details, see \cite{OT}).
All sets in Table 1.1 except $S_{13}$ was firstly provided by Jagy and he checked that each quadratic form in the sets represents same integers up to $1,000,000$.
\begin{table}[h]
\centering
\renewcommand{\arraystretch}{1}\renewcommand{\tabcolsep}{1mm}
\begin{tabular}{l|ll}\multicolumn{3}{c}{\textbf{Table 1.1.} Ternary quadratic forms with same representations}\\
\hline
$S_1$ &  $5x^2+8y^2+8z^2-5yz-xz-4xy$, & $5x^2+5y^2+8z^2-yz-4xz-2xy$ \\ \hline
$S_2$ & $3x^2+4y^2+7z^2-yz$, & $3x^2+4y^2+7z^2+4yz+3xz+3xy$  \\ \hline
$S_3$ & $x^2+4y^2+7z^2-yz$, & $x^2+4y^2+5z^2-yz-xz$  \\ \hline
$S_4$ &  $4x^2+7y^2+25z^2-4yz-2xz-2xy$, & $4x^2+7y^2+7z^2+5yz+2xz+2xy$ \\ \hline
$S_5$ & $2x^2+6y^2+41z^2-3yz-xz$, & $2x^2+2y^2+41z^2+yz+2xz+2xy$  \\ \hline
$S_6$ &  $2x^2+6y^2+14z^2-3yz-xz$, & $2x^2+2y^2+14z^2+yz+2xz+2xy$ \\ \hline
$S_7$ & $2x^2+4y^2+8z^2+4yz+xz+xy$, & $2x^2+2y^2+4z^2-yz-2xz$  \\ \hline
$S_8$ & $5x^2+5y^2+8z^2-4xz-3xy$, & $5x^2+7y^2+7z^2+6yz+xz+5xy$  \\ \hline
$S_9$ &  $3x^2+3y^2+7z^2+yz+2xz+xy$, & $3x^2+5y^2+5z^2+3yz+xz+3xy$ \\ \hline
$S_{10}$ & $5x^2+5y^2+8z^2-yz-2xz-4xy$, & $5x^2+5y^2+6z^2-3xz-2xy$  \\ \hline
$S_{11}$ & $2x^2+4y^2+7z^2-xz-xy$, & $2x^2+4y^2+7z^2+4yz+2xz+xy$  \\ \hline
$S_{12}$ & $4x^2+6y^2+7z^2+3yz+2xz+3xy$, & $4x^2+4y^2+6z^2-3xz-2xy$  \\ \hline
$S_{13}$ & $5x^2+12y^2+28z^2-4xz-4xy$, & $5x^2+12y^2+24z^2-8yz-4xy$, \\  
~ & $5x^2+12y^2+21z^2-4yz-2xz-4xy$, & $5x^2+12y^2+12z^2-4xz-4xy$ \\ \hline
$S_{14}$ & $3x^2+5y^2+7z^2-2yz-2xy$, & $3x^2+5y^2+6z^2-2xz-2xy$, \\
~ & $3x^2+5y^2+6z^2+4yz+2xz+2xy$, & $3x^2+3y^2+5z^2-2yz-2xz$ \\ \hline
$S_{15}$ & $3x^2+5y^2+5z^2+5yz+2xz+3xy$, & $3x^2+3y^2+5z^2-yz-2xz-xy$,\\
~& $3x^2+3y^2+5z^2+2yz+3xz+xy$, & $3x^2+3y^2+3z^2+yz+xz+3xy$\\ \hline
\end{tabular}
\end{table}

In this article, we prove the existence of pairs of ternary quadratic forms representing same integers which are not in the Kaplansky's list. Actually, we prove that each of two quadratic forms contained in $S_4$, $S_6$, $S_7$ and $S_8$ in Table 1.1 represents same integers.

Let 
$$
f(x_1,x_2,\dots,x_r)=\sum_{1 \le i, j\le r} a_{ij} x_ix_j \ (a_{ij}=a_{ji} \in \mathbb Z)
$$
 be a positive definite integral quadratic form. 
The corresponding integral symmetric matrix of $f$ is defined by $M_f=(a_{ij})$. 
For an integer $k$, we say $k$ is {\it represented by $f$} if the equation $f(x_1,x_2,\dots,x_m)=k$ has an integer solution. 

Any unexplained notations and terminologies can be found in \cite{ki1} or \cite{om}.

\section{General tools}

In \cite{4-8}, \cite{regular}, and \cite{pentagonal}, we developed a method that determines whether or not integers in an arithmetic progression are represented by some particular ternary quadratic form. 
We briefly introduce this method for those who are unfamiliar with it. 
 
 Let $d$ be a positive integer and let $a$ be a nonnegative integer $(a \leq d)$. 
We define 
$$
S_{d,a}=\{dn+a \mid n \in \mathbb N \cup \{0\}\}.
$$
For two integral ternary quadratic forms $f,g$, we define
$$
R(g,d,a)=\{v \in (\mathbb{Z}/d\mathbb{Z})^3 \mid vM_gv^t\equiv a \ (\text{mod }d) \}
$$
and
$$
R(f,g,d)=\{T\in M_3(\mathbb{Z}) \mid  T^tM_fT=d^2M_g \}.
$$
A coset (or, a vector in the coset) $v \in R(g,d,a)$ is said to be {\it good} with respect to $f,g,d, \text{ and }a$ if there is a $T\in R(f,g,d)$ such that $\frac1d \cdot vT^t \in \mathbb{Z}^3$.  
The set of all good vectors in $R(g,d,a)$ is denoted by $R_f(g,d,a)$.   
If  $R(g,d,a)\setminus R_f(g,d,a)=\emptyset$, then we write  
$$
g\prec_{d,a} f.
$$ 
If $g\prec_{d,a} f$, then by Lemma 2.2 of \cite{regular},  we have 
\begin{equation}\label{good}
S_{d,a}\cap Q(g) \subset Q(f).
\end{equation}

In general, if $d$ is large, then it is hard to compute the set  $R(g,d,a)\setminus R_f(g,d,a)$ exactly by hand. 
A MAGMA based  computer program for computing the set could be available upon request to the author.

\section{Representations of ternary quadratic forms}

In this section, we prove that each of two quadratic forms contained in $S_4$, $S_6$, $S_7$ and $S_8$ in Table 1.1 represents same integers. 

For some technical reason, we consider integral quadratic forms obtained from quadratic forms in Table 1.1 by scaling 2. 
Let $m$ be a positive integer.
Note that for two quadratic forms $f$ and $g$,  $Q(f)=Q(g)$ if and only if  $Q(mf)=Q(mg)$.
\vskip0.5pc
\noindent\textbf{(i)} Quadratic forms in $S_4$. Let $f(x,y,z)=8x^2+14y^2+50z^2-8yz-4xz-4xy$ and  $g(x,y,z)=8x^2+14y^2+14z^2+10yz+4xz+4xy$.
Then each symmetric matrix corresponding to $f$ and $g$ is as follows: 
$$
M_f=\begin{pmatrix}8&-2&-2\\-2&14&-4\\-2&-4&50\end{pmatrix}\quad\text{and}\quad M_g=\begin{pmatrix}8&2&2\\2&14&5\\2&5&14\end{pmatrix}.
$$
One may easily show that for $T=\begin{pmatrix}1&0&0\\0&0&-2\\0&-1&1 \end{pmatrix}$, $T^tM_gT=M_f$, which implies that $f$ is a subform of $g$. 
Therefore, $Q(f)\subset Q(g)$. 

It is enough to show that every positive even integer represented by $g$ is also represented by $f$.
One may easily compute that 
$$
R(g,4,0)=\{(v_1,v_2,v_3)\in(\mathbb{Z}/4\mathbb{Z})^3 \mid v_2\equiv0 \Mod2, v_3\equiv0\Mod2\}.
$$
Furthermore, note that $R(f,g,4)$ consists of $8$ matrices, and in particular, it contains the following matrix:
$$
T_1=\begin{pmatrix}4&2&2\\0&4&2\\0&0&2\end{pmatrix}.
$$
Then, we know that $\frac14\cdot(v_1,v_2,v_3)\cdot T_1^t\in\mathbb{Z}^3$, for any $(v_1,v_2,v_3)\in R(g,4,0)$. 
Therefore, $g\prec_{4,0}f$ holds. 
Then by \eqref{good}, every positive integer congruent to $0$ modulo $4$ which is represented by $g$ is also represented by 
$f$.
As a sample, for $(1,4,2)\in\mathbb{Z}^3$, note that $g(1,4,2)=392$ and $(1,4,2)$ is congruent to $(1,0,2)$ modulo $4$ which is contained in $R(g,4,0)$. 
Then
$$
f(4,5,1)=f\left(\frac14\cdot(1,4,2)\cdot T_1^t\right)=g(1,4,2)=392.
$$
Similarly, one may easily show that 
$$
g\prec_{12,6}f,\quad g\prec_{12,10}f.
$$
In fact, $Q(g)\cap\{12n+10\mid n\in\mathbb{N}\cup\{0\}\}=\emptyset$.
However, we consider this case for completeness of the proof. 
It also happens for some other cases that we will consider later.
Therefore, by \eqref{good} every positive integer congruent to $6$ or $10$ modulo $12$ which is represented by $g$ is also represented by $f$.
Now, for some nonnegative integer $n_1$ and $v\in\mathbb{Z}^3$, assume that $g(v)=12n_1+2$. 
One may easily check that there are exactly $864$ vectors in $R(g,12,2)$ and $144$ matrices in $R(f,g,12)$. 
Furthermore, all vectors in $R(g,12,2)$ are good vectors with respect to $f,g,12$ and $2$ except $32$ vectors. 
Note that 
$$
\begin{array}{l}
\!\!\!R(g,12,2)\setminus R_f(g,12,2)\\
=\{(v_1,v_2,v_3)\in(\mathbb{Z}/12\mathbb{Z})^3 \mid v_1\not\equiv0\Mod3, v_2\equiv\pm3\Mod{12}, v_3\equiv\pm3\Mod{12}\}.
\end{array}
$$
Now, define 
$$
\widetilde{T}=\begin{pmatrix}12&6&2\\0&0&12\\0&-12&-8\end{pmatrix}.
$$
Note that $\widetilde{T}^tM_g\widetilde{T}=12^2M_g$.
Then one may easily check that for any vector $u\in \mathbb{Z}^3$ such that $u\Mod{12}\in R(g,12,2)\setminus R_f(g,12,2)$,
\begin{equation}\label{integer vector}
\frac{1}{12}\cdot u\widetilde{T}^t\in\mathbb{Z}^3.
\end{equation}
All computations were done by a computer program based on MAGMA. 
If $v\Mod{12}$ is a good vector with respect to $f,g,12$ and $2$, then there is a $T_2\in R(f,g,12)$ such that $\frac{1}{12}\cdot vT_2^t\in\mathbb{Z}^3$ and $f(\frac{1}{12}\cdot vT_2^t)=12n_1+2$. 
Assume that $v\Mod{12}$ is not a good vector.
Then by \eqref{integer vector}, we know that $\frac{1}{12}\cdot v\widetilde{T}^t\in\mathbb{Z}^3$. 
If $\frac{1}{12}\cdot v\widetilde{T}^t\Mod{12}\in R_f(g,12,2)$, then we are done. 
If $\frac{1}{12}\cdot v\widetilde{T}^t\Mod{12}$ is not good, then $\left(\frac{1}{12}\right)^2\cdot v(\widetilde{T}^t)^2\in\mathbb{Z}^3$ by \eqref{integer vector}.
Now inductively, we may assume that for any positive integer $m$, 
$$
\left(\frac{1}{12}\right)^m \cdot v (\widetilde{T}^t)^m\in \mathbb{Z}^3 \text{ and } \left(\frac{1}{12}\right)^m \cdot v (\widetilde{T}^t)^m \Mod{12}\in R(g,12,2)\setminus R_f(g,12,2).
$$ 
Since there are only finitely many integer solution of $g(x,y,z)=12n_1+2$ and $\frac{1}{12}\widetilde{T}$ has an infinite order, it is impossible, unless $v$ is an eigenvector of $\widetilde{T}$. 
Note that $(\pm1,0,0)$ are the only integral primitive eigenvectors of $\widetilde{T}$. 
Hence, if $12n_1+2$ is not of the form $g(\pm t,0,0)=8t^2$ for some positive integer $t$, then it is also represented by $f$. 
Furthermore, $8t^2$ is represented by $f$ since $8$ is represented by $f$. 
Therefore $12n_1+2$ is represented by $f$. 
Hence, every positive even integer which is represented by $g$ is also represented by $f$.
Therefore, $f$ and $g$ have perfectly identical integral representations.
\vskip0.5pc
\noindent\textbf{(ii)} Quadratic forms in $S_6$. Let $f(x,y,z)=4x^2+12y^2+28z^2-6yz-2xz$ and $g(x,y,z)=4x^2+4y^2+28z^2+2yz+4xz+4xy$.
One may easily show that $f$ is a subform of $g$. 
Therefore $Q(f)\subset Q(g)$.

It is enough to show that every positive even integer represented by $g$ is also represented by $f$. 
One may easily show that 
$$
g\prec_{4,2}f,\quad g\prec_{8,0}f, \quad g\prec_{24,12}f, \quad g\prec_{24,20}f \quad g\prec_{48,4}f, \quad\text{and}\quad g\prec_{48,28}f.
$$
Therefore, we know that $Q(g)\subset Q(f)$ by \eqref{good}.
Hence, $f$ and $g$ have perfectly identical integral representations.
\vskip0.5pc
\noindent\textbf{(iii)} Quadratic forms in $S_7$. Let $f(x,y,z)=4x^2+8y^2+16z^2+8yz+2xz+2xy$ and $g(x,y,z)=4x^2+4y^2+8z^2-2yz-4xz$.
One may easily show that $f$ is a subform of $g$. Therefore, $Q(f)\subset Q(g)$.

It is enough to show that every positive even integer represented by $g$ is also represented by $f$. 
One may easily show that 
$$
\begin{array}{llll}
g\prec_{4,2}f, & g\prec_{24,0}f, & g\prec_{24,4}f, &~\\
g\prec_{24,8}f, & g\prec_{24,12}f, & g\prec_{24,16}f, & \text{and}\quad g\prec_{24,20}f.
\end{array}
$$
Therefore, we know that $Q(g)\subset Q(f)$ by \eqref{good}.
Hence, $f$ and $g$ have perfectly identical integral representations.
\vskip0.5pc
\noindent\textbf{(iv)} Quadratic forms in $S_8$. Let  $f(x,y,z)=10x^2+10y^2+16z^2-8xz-6xy$ and $g(x,y,z)=10x^2+14y^2+14z^2+12yz+2xz+10xy$ .
One may easily show that
$$
g\prec_{4,0}f, \quad g\prec_{12,2}f, \quad g\prec_{12,6}f, \quad
g\prec_{36,10}f, \quad g\prec_{36,22}f, \quad\text{and}\quad g\prec_{36,34}f.
$$
Therefore, we know that $Q(g)\subset Q(f)$ by \eqref{good}.

Conversely, one may easily show that
$$
f\prec_{4,0}g, \quad f\prec_{12,2}g, \quad f\prec_{12,6}g,\quad
f\prec_{36,10}g, \quad f\prec_{36,22}g, \quad\text{and}\quad f\prec_{36,34}g.
$$
Therefore, we know that $Q(f)\subset Q(g)$ by \eqref{good}.
Hence, $f$ and $g$ have perfectly identical integral representations.

\begin{rmk}
{\rm(1) For the other pairs in Table 1.1 not considered in this article, it seems to be inadequate to apply our method.}

\noindent{\rm(2) In 1981, Hsia conjectured in \cite{hsia} that two positive ternary quadratic forms in the same genus are equivalent if they primitively represent the same integers. 
Jagy found two candidates of pairs of ternary quadratic forms violating Hsia's conjecture.
Actually, He proved that each of two quadratic forms contained in $S_8$ and $S_9$ in Table 1.1 which are contained in same genus primitively represents same integers up to $1,000,000$ (see \cite{jagy}).} 
\end{rmk}

\end{document}